\documentclass[12pt]{amsart}

\usepackage{amssymb}
\usepackage{amsmath}
\evensidemargin\oddsidemargin

\begin{document}


\setcounter{page}{1}

\newtheorem{Eins}{Lemma 1\!\!}
\newtheorem{Zwei}{Remark 2\!\!}
\newtheorem{Drei}{Lemma 3\!\!}
\newtheorem{Vier}{Remarks 4\!\!}
\newtheorem{Fuenf}{Remarks 5\!\!}
\newtheorem{Korse}{Corollary 6\!\!}
\newtheorem{Prosi}{Proposition 7\!\!}
\newtheorem{Acht}{Final remarks 8\!\!}
\newtheorem{THEO}{Theorem\!\!}

\renewcommand{\theEins}{}
\renewcommand{\theZwei}{}
\renewcommand{\theDrei}{}
\renewcommand{\theVier}{}
\renewcommand{\theFuenf}{}
\renewcommand{\theKorse}{}
\renewcommand{\theProsi}{}
\renewcommand{\theAcht}{}
\renewcommand{\theTHEO}{}

\def\a{\alpha}
\def\b{\beta}
\def\CC{{\mathbb{C}}} 
\def\cV{{\mathcal V}}
\def\Ea{E_\a}
\def\EE{{\mathbb{E}}} 
\def\eps{\varepsilon}
\def\fa{f_\a}
\def\far{f_\a^r}
\def\Fa{F_\a}
\def\Far{F_{\a,r}}
\def\Fat{G_{\a,t}}
\def\fda{{\tilde f}_\a}
\def\ga{g_\a}
\def\gar{g_\a^r}
\def\Ga{G_\a}
\def\har{h_\a^r}
\def\Hta{H_{\a,t}}
\def\ii{{\rm i}}
\def\lbd{\lambda}
\def\lacc{\left\{}
\def\lcr{\left[}
\def\lpa{\left(}
\def\lva{\left|}
\def\LL{\mathcal L}
\def\Ma{M_\a}
\def\MM{\mathcal M}
\def\NN{{\mathbb{N}}} 
\def\pb{{\mathbb{P}}}
\def\rl{{\mathbb{R}}}
\def\racc{\right\}}
\def\rcr{\right]}
\def\rpa{\right)}
\def\Rde{\tilde R}
\def\Rta{\hat R}
\def\RR{\mathcal R}
\def\Tde{\tilde T}
\def\Tta{\hat T}
\def\Va{V_\a}
\def\Xar{X_{\a,r}}
\def\Yar{Y_{\a,r}}
\def\Ya{Y_\a}
\def\Zba{{\bar Z}_\a}
\def\Zda{{\tilde Z}_\a}
\def\Xda{{\tilde X}_\a}
\def\Za{Z_\a}
\def\Zar{Z_\a^r}
\def\Ua{U_\a}
\def\UU{\mathcal U}

\def\d{\, \mathrm{d}}
\def\elaw{\stackrel{d}{=}}

\newcommand{\fin}{\vspace{-0.3cm}
                  \begin{flushright}
                  \mbox{$\Box$}
                  \end{flushright}
                  \noindent}
              
\title[Unimodality and stable densities]{On the unimodality of power transformations of positive stable densities}

\author[Thomas Simon]{Thomas Simon}

\address{Laboratoire Paul Painlev\'e, U. F. R. de Math\'ematiques, Universit\'e de Lille 1, F-59655 Villeneuve d'Ascq Cedex. {\em Email} : {\tt simon@math.univ-lille1.fr}}



\begin{abstract}  Let $\Za$ be a positive $\a-$stable random variable and $r\in\rl.$ We show the existence of an unbounded open domain $D$ in $[1/2,1]\times\rl$ with a cusp at $(1/2,-1/2)$, characterized by the complete monotonicity of the function $\Far (\lbd) = (\a \lbd^\a -r)e^{-\lbd^\a}\!\! ,$ such that $\Za^r$ is unimodal if and only if $(\a, r)\notin D.$  
\end{abstract}

\maketitle


\vspace{5mm}

\section{Introduction and statement of the result}

A real random variable $X$ is said to be unimodal (or quasi-concave) if there exists $a\in\rl$ such that its distribution function $\pb[X\le x]$ is convex on $(-\infty, a)$ and concave on $(a, +\infty)$. The law of $X$ decomposes then into $\pb[X\in dx] = c \delta_a (dx) + f(x) dx,$ where $f$ is non-decreasing on $(-\infty, a)$ and non-increasing on $(a, +\infty).$ The number $a$ is called a mode of $X$ and is not necessarily unique. We will denote by $\UU_a$ the set of unimodal random variables with a mode at $a$ and set $\UU = \cup_{a\in\rl}\;\UU_a.$ The set of random variables with a monotone density function, which is included in $\UU,$ will be denoted by $\MM.$ 

Unimodal random variables share many interesting properties in mathematical statistics, sometimes analogous to those of the normal distribution. A celebrated example is the following three-rule: when $X$ is unimodal and square-integrable, then
$$\pb[\vert X -x \vert > 3 \tau_x]\;\le\; 4/81\; <\; 0.05$$ 
for all $x\in\rl,$ with the notation $\tau_x^2 = \EE[(X-x)^2].$ This optimal bound was discovered by Gauss in 1823 when $x$ is a mode of $X,$ and extended to all $x$ by P\'etunine \& Vissotchanski in 1983.  Notice that in the particular case $x= \EE[X]$ the above rule divides the classical bound of Bienaym\'e \& Tchebitcheff in more than one-half. We refer to the monograph \cite{DJP} for details, references and further results concerning unimodality.

In the present paper the unimodality of $Z_\a^r$ is discussed, where $\Za$ is a positive $\a-$stable random variable ($0<\a<1$) with density function $\fa$ normalized such that
\begin{equation}
\label{normal}
\int_0^\infty e^{-\lbd t} \fa(t) dt\; =\; \EE\lcr e^{-\lbd \Za}\rcr \; =\; e^{-\lbd^\a}, \quad \lbd \ge 0,
\end{equation}
and $r\in\rl.$ Apart from the trivial situation $r=0$ the random variable $\Zar$ has always a smooth density $\far$ and we would like to know if $\far$ has a unique local maximum value, or not. There are already some results in the literature, all with a positive answer.

\vspace{2mm}

$\bullet$ First of all, the cases $r =1$ and $\{\a \in [1/2,1), r = -\a\}$ were considered in \cite{CI}. When $r=1$ the fact that $\Za\in\UU$ is a consequence of Lemma 1 and the proof of the Theorem therein. The other case comes after recalling - see e.g. Exercises 29.7 and 29.18 in \cite{S}, that 
$$Z_\a^{-\a}\elaw X_\a\; \mbox{conditioned on} \; X_\a >0,$$ 
where $X_\a$ is a $(1/\a)-$stable random variable with skewness parameter $\beta = -1,$ so that one can apply Lemma 1' and a similar discussion to the proof of the Theorem  in \cite{CI} to get $Z_\a^{-\a}\in\UU.$

\vspace{2mm}

$\bullet$ The case $r \le \a/(\a-1)$ is a corollary to the main theorem in \cite{PSS}, which shows a stronger property. From (3.1) therein one obtains indeed the following representation
$$\far(x)\; =\; \int_0^\infty e^{-xy} y\, \pb[\exp\{V_{\a, -r}\}\in dy],$$
where the positive random variable $V_{\a, -r}$ has Laplace transform given by
$$\EE [e^{-s V_{\a, -r}}] \; =\; \frac{\Gamma(1+sr/\a)}{\Gamma(1+s)\Gamma(1+sr)}$$
and is actually infinitely divisible (ID). The easy part of Bernstein's theorem - see e.g. Theorem 1.4. in \cite{SSV} - entails then that $\far$ is completely monotone (CM), hence monotone, and one gets $\Zar\in\UU_0.$ Notice that $\Zar$ is ID as well - see Theorem 51.6 in \cite{S}.

\vspace{2mm}

$\bullet$ Last, the case $r =\a$ follows from the remark made in \cite{Pat} - see Paragraph 3.2 therein, that $Z_\a^\a$ is self-decomposable (SD): the classical result of Yamazato - see e.g. Theorem 59.12 in \cite{S} - yields then $Z_\a^\a\,\in\,\UU.$

\vspace{2mm}

However, one cannot deduce further positive results from the above three cases, because unimodality is hardly stable under power transformations. We may also recall that proving or disproving unimodality can actually turn out to be a difficult problem for random variables whose density function cannot be written in closed form. Our main result answers the question for $\Zar$ in the following way:

\begin{THEO} {\rm (a)} One has $\Zar\in \UU\cap\MM^c$ for all $r > -\a.$ 

\vspace{2mm}

{\rm (b)} One has $Z_\a^{-\a}\in \MM$ for $\a\le 1/2$ and $Z_\a^{-\a}\in \UU\cap\MM^c$ for $\a > 1/2.$

\vspace{2mm}

{\rm (c)} In the case $r < -\a$ there exists a homeomorphism $R : [1/2,1) \to [1/2, +\infty)$ such that
$$\Zar\,\in\,\UU\;\Leftrightarrow\;\Zar\,\in\,\MM\;\Leftrightarrow\;\Far\;\, \mbox{\em is CM}\;\Leftrightarrow\;\a\le 1/2\;\, {\rm or}\;\, r \le -R(\a),$$
with the notation $\Far (\lbd) = (\a \lbd^\a -r)e^{-\lbd^\a}\!\!$ for all $\lbd \ge 0.$ Moreover, one has the following bounds:
$$1/4(1-\a)\;\le\; R(\a)\;\le\; (\a/\sin^2(\pi \a))\wedge (\a/(1-\a)),$$
so that $R(\a)\sim\a$ as $\a \to 1/2$ and $R(\a)\asymp 1/(1-\a)$ as $\a \to 1.$ In particular, the open domain $D = \{(\a,r)\; /\; \Zar\notin\UU\}$ is unbounded and has a cusp at $(1/2,-1/2).$
\vspace{1mm}

\end{THEO}

One can understand from the limiting case $r = -\a$ why there should be a prohibited zone for the unimodality of $\Zar$ in the area $\{\a >1/2, r < -\a\}.$ The well-known Humbert-Pollard series representation for $\fa$ - see (\ref{HP}) below - yields indeed the following limiting values
$$f_\a^{-\a}(0+)\; =\; \frac{1}{\Gamma(1-\a)}\; , \;\; (f_\a^{-\a})'(0+)\; =\; \frac{-2}{\Gamma(1-2\a)}\; \;\mbox{and} \;\; f_\a^r(0+)\; =\; +\infty\;\mbox{for all $r < -\a.$}$$ 
In particular, the random variable $Z_\a^{-\a}$ has a mode which is not zero for $\a >1/2,$ and by the continuity of $r\mapsto \far(x)$ for all $x >0$ one sees that $\far$ must have at least two separate modes (one of which is zero) provided that $r$ is close enough to $-\a.$ However, it is somewhat surprising for the forbidden area $D$ to be cuspidal at $(1/2, -1/2)$ and we could not find any quick, intuitive explanation for that. As in \cite{Th}, it also appears that the cases $\a\le 1/2$ and $\a > 1/2$ are completely different as far as unimodality properties for $\Za$ are concerned. When $\a \le 1/2,$ the condition $r\le -\a$ characterizes namely the complete monotonicity of $\Far,$ but also its weaker log-convexity as well as the much stronger property that  $-\log \Far$ is a Bernstein function (in other words, that $\Far$ is the Laplace transform of a positive ID law). When $\a > 1/2$ these three equivalent characterizations do not hold anymore, and finding a closed expression in $\a$ for the frontier function $R(\a)$ seems somehow hopeless. 

From the limiting case $\{\a >1/2, r = -\a\}$ it is reasonable to conjecture that the variable $\Zar$ is bimodal when $(\a, r) \in D.$ In general, techniques for investigating multimodality are different from those devoted to unimodality. Anyway it does not seem that the methods of this paper, which in the non-monotonic situation rely on the notion of strong multiplicative unimodality \cite{CT}, can be of any help for proving bimodality. It would also be interesting to study the unimodality of $X_\a^r,$ where $X_\a$ is a general $\a-$stable variable conditioned to stay positive. Such cut-off variables have been introduced in Chapter 3 of \cite{Z} in the framework of $M$-infinite divisibility. Their multiplicative factorizations are however more complicated than in the true positive situation, and we will hence leave this question open to some further research, as well as that of the bimodality of $\Zar$ when $(\a,r)\in D.$

\section{Proof of the theorem}

\subsection{Two lemmas} For parts (a) und (b) we will use the following multiplicative factorization of $\Za$ which was discovered by Kanter - see Corollary 4.1. in \cite{K} - as a direct consequence of a contour integration made by Chernine \& Ibragimoff - see the final remark in \cite{CI}: one has
\begin{equation}
\label{Kant}
\Za\; \elaw\; L^{(\a-1)/\a}\times b_\a^{-1/\a}(U),
\end{equation}
where $L\sim {\rm Exp} (1),$ $U\sim {\rm Unif}(0,\pi)$ independent of $L$, and 
$$b_\a (u) = (\sin u/\sin (\a u))^\a(\sin u/\sin ((1-\a) u))^{1-\a}, \quad u\in(0,\pi).$$
We will need the following property of $b_\a$, partly already shown in \cite{K}.

\begin{Eins} The function $b_\a$ is decreasing and concave on $(0,\pi).$
\end{Eins}

\begin{proof}
First, it follows directly from the beginning of the proof of Theorem 4.1. in \cite{K} that $b_\a$ decreases and is log-concave on $(0,\pi).$ When $\a =1/2,$ one can write $b_{1/2}(u) = 2\cos(u/2)$ which is surely a concave function on $(0,\pi).$  When $\a\neq 1/2,$ it is however more difficult to prove that $b_\a$ is concave. We will suppose $0<\a<1/2$ and set $\b = 1-\a.$ A computation yields 
$$b_\a''(u)\; =\; b_\a(u) (\a(A_\a^2(u) -A_\a'(u)) +\b (A_\b^2(u) -A_\b'(u)) -\a\b (A_\a(u) - A_b(u))^2),$$
where $A_\a (u) = \a \cot (\a u) - \cot (u)$ and $A_\b (u) = \b \cot (\b u) - \cot (u).$ It is hence sufficient to show that the quantity 
$$\a(A_\a^2(u) -A_\a'(u)) +\b (A_\b^2(u) -A_\b'(u))\; = \; 
2(\a^2\cot(\a u) A_\a (u) +\b^2\cot (\b u) A_\b(u))- 3\a\b$$
is non-positive for all $u\in (0, \pi).$ We will obtain the slightly better inequality
\begin{eqnarray}
\label{Haupt}
\a^2\cot(\a u) A_\a (u) +\b^2\cot (\b u) A_\b(u) & \le & \a\b
\end{eqnarray}
and for this purpose we change the expression on the left-hand side into 
$$\a\b (\cot(\a u) A_\a (u) +\cot (\b u) A_\b(u)) \; +\; (\b -\a) (\b\cot (\b u) A_\b(u)  -\a\cot(\a u) A_\a (u)).$$
A further decomposition gives
\begin{eqnarray*}
\cot(\a u) A_\a (u) +\cot (\b u) A_\b(u) & = & \a\cot^2(\a u) +\b\cot^2 (\b u) A_\b(u) + 1 -\cot(\a u)\cot (\b u)\\
& = & (A_\a (u) - A_\b (u))(\cot(\a u) -\cot (\b u)) + 1
\end{eqnarray*}
and
\begin{eqnarray*}
\b\cot (\b u) A_\b(u)  -\a\cot(\a u) A_\a (u) & = & (A_\a (u) - A_\b (u))(\cot (u) - \cot(\a u) -\cot (\b u)).
\end{eqnarray*}
After some simplification one finds that (\ref{Haupt}) is equivalent to
\begin{equation}
\label{Neben}
(A_\a (u) - A_\b (u))(\a A_\a (u) - \b A_\b (u))\; \le \; 0,\quad  u\in(0,\pi).
\end{equation}
Differentiating in $\a$ it is easy to see - and already used in the proof of Theorem 4.1. in \cite{K} - that $A_\a (u) \ge A_\b(u)$ for every $u\in (0,\pi).$ It is however more painful to show $\a A_\a \le \b A_\b$ on $(0,\pi)$ with the help of successive derivatives, and we will rather appeal to the eulerian formula
$$\pi \cot(\pi z) \; =\; \frac{1}{z} \; +\; 2z \sum_{n\ge 1} \frac{1}{z^2 - n^2}\cdot$$  
Setting $u =\pi z,$ the latter entails namely
\begin{eqnarray*}
\a A_\a (u) - \b A_\b (u) & = & \frac{2\a\b (\b^2 - \a^2)z^3}{\pi}\sum_{n\ge 1}\frac{n^2}{(z^2 -n^2)(\a^2z^2 -n^2)(\b^2 z^2 -n^2)} \; \le \; 0,
\end{eqnarray*}
the inequality being justified by $\b \ge \a$ and $z\in (0,1).$ This shows (\ref{Neben}) and completes the proof.

\end{proof}

\begin{Zwei} {\em It is easy to see that $G_{\a,\b} (u) = \a^2\cot(\a u) A_\a (u) +\b^2\cot (\b u) A_\b(u) \to \a\b$ as $u\to 0,$ so that the lemma would follow as soon as it is shown that $G_{\a, \b}$ is non-increasing on $(0,\pi).$ This latter property can be obtained similarly as above but we could not find any simpler method.}
\end{Zwei}

For part (c) we will need the following interesting property of $\Far,$ whose proof relies partly on the theory of Bernstein functions - see \cite{SSV} for a modern account, especially Chapter 3 therein.

\begin{Drei} There exists an increasing function $R : [0,1] \to [0, +\infty],$ such that
\begin{equation}
\label{WM}
\Far\;\, \mbox{\em is CM}\;\Longleftrightarrow\; r \le -R(\a).
\end{equation}
Moereover, one has $R(\a) = \a$ for all $\a \in [0, 1/2]$ and the bounds $1/4(1-\a)\,\le\, R(\a)\,\le\, \a/\sin^2(\pi \a)$ hold for all $\a \in [1/2, 1].$
\end{Drei}

\begin{proof} Let us first consider the function $\Fat = (\lbd^\a + t)e^{-\lbd^\a}$ on $\rl^+.$  Since $e^{-\lbd^\a}$ itself is CM - recall (\ref{normal}) - we see that the notation $T(\a) = \inf\{ t\in \rl, \; \Fat\,\mbox{is CM}\}$ makes sense and - from the continuity in $t$ of the successive derivatives of $\Fat$, that
$$\Fat\;\mbox{is CM}\;\Longleftrightarrow\; t \ge T(\a).$$ Moreover it is clear that $T(0) = 0, T(\a) \ge 0$ and from the computation 
$$G_{1,t}^{(n)} (\lbd) = (-1)^n e^{-\lbd} (\lbd + t -n),$$ 
that $T(1) = +\infty.$ Next we show indirectly that $T$ is non-decreasing on $[0,1]$: supposing $\a_1 < \a_2$ and $T(\a_1) > T(\a_2)$ and considering then $t \in (T(\a_2), T(\a_1))$ and $\gamma = \a_1/\a_2\in [0,1),$ because $G_{\a_2, t}$ is CM and $\lbd^\gamma$ is a  Bernstein function - see again (\ref{normal}), it follows from Criterion 2 p. 417 in \cite{F} - see also Theorem 3.6. in \cite{SSV} - that $G_{\a_1, t}$ is also CM, which contradicts the definition of $T(\a_1),$ so that $T$ does not decrease on $[0,1]$. Finally, the function $R(\a) = \a T(\a)$ increases on $[0,1],$ fulfils $R(0) = 0, R(1) = +\infty,$ and it is clear from the above discussion that the equivalence (\ref{WM}) holds.

For the lower bound we first recall that by Bernstein's theorem and H\"older's inequality, the CM property implies the log-convexity one. A computation yields on the other hand
 $$\Fat\;\mbox{is log-convex on $\rl^+$}\;\Longleftrightarrow\; x^2 + (2t - 1/(1-\a))x + t^2 -t \ge 0\;\mbox{for all $x\ge 0.$}$$ 
This entails that if $\Fat$ is CM, then necessarily $t^2\ge t$ and either $t \ge 1/2(1-\a),$ or $(2t - 1/(1-\a))^2 \le 4(t^2-t),$ from which the bounds $R(\a)\ge \a$ for all $\a\in [0,1/2]$ and $R(\a)\ge 1/4(1-\a)$ for all $\a\in [1/2,1]$ easily follow.

For the upper bound we finally notice that if $\Hta(\lbd) =\lbd^\a - \log (1+\lbd^\a/t)$ is a Bernstein function, then $\Fat$ is CM as $t^{-1}$ times the Laplace-transform of an ID law. Besides, $\lbd^\a$ and $\log (1+\lbd^\a/t)$ are Bernstein functions themselves, whose L\'evy measures are given by
$$\lbd^\a\; =\; \frac{\a}{\Gamma(1-\a)}\int_0^\infty\!\!\! (1- e^{-\lbd x})\frac{dx}{x^{1+\a}}\;\;\mbox{and}\;\;\log (1 + \lbd^\a/t)\; =\; \a\int_0^\infty\!\!\! (1- e^{-\lbd x})\Ea(-t x^\a)\frac{dx}{x}\cdot$$
Above,
$$\Ea(x)\; =\; \sum_{n=0}^{\infty} \frac{x^n}{\Gamma (1 +\a n)}$$
is the Mittag-Leffler function with index $\a$ and the second computation comes from Remark 2.2. in \cite{P}, having made the correction $x^k \to x^{\a k}$ therein. Recall also the asymptotic behaviour
$$\Ea(-x^\a) \;\sim \; \frac{x^{-\a}}{\Gamma(1-\a)}$$
as $x\to +\infty$ - see (18.1.7) in \cite{E}, which ensures the convergence of the integral in the second formula. One gets
\begin{equation}
\label{MrL}
\Hta(\lbd)\; =\; \frac{\a}{t\Gamma(1-\a)}\int_{\rl^+} (1- e^{-\lbd x})(t- \Ua(t x^\a))\frac{dx}{x^{1+\a}}
\end{equation}
with the notation $\Ua(x) = \Gamma(1-\a) x \Ea(-x),$ so that $\Hta$ is Bernstein if and only if
\begin{equation}
\label{You}
\Ua (x^\a) \le t
\end{equation}  
for all $x\ge 0.$ From Exercise 29.18 in \cite{S}, (\ref{normal}) in the present paper and Exercise 4.21 (3) in \cite{CY} one obtains the classical integral representation 
\begin{equation}
\label{ML}
\Ea(-x^\a)\; =\; \EE \lcr e^{-x^\a/\Za^\a}\rcr \; = \; \EE\lcr e^{-x (\Zba/\Za)}\rcr\; =\; \frac{\sin \pi \a}{\pi}\int_{\rl^+} \frac{u^{\a-1} e^{-xu} }{u^{2\a} + 2u^\a\cos \pi\a +1} \, du,
\end{equation}
where in the second equality $\Zba$ is an independent copy of $\Za,$ and the complement formula for the Gamma function entails 
$$\Ua(x^\a) \; = \; \frac{1}{\Gamma(\a)}\int_0^\infty \!\frac{u^{\a-1} e^{-u} }{(u/x)^{2\a} + 2(u/x)^\a\cos \pi\a +1} \, du,$$
an expression which is everywhere smaller than 1 when $\a\le 1/2$ and everywhere smaller than $1/\sin^2(\pi\a)$ when $\a \ge 1/2.$ In other words, one has $T(\a) \le 1$ when $\a\le 1/2$ and $T(\a) \le 1/\sin^2(\pi\a)$ when $\a \ge 1/2.$ 

These two bounds yield finally $R(\a) =\a$ for all $\a\in [0,1/2]$ and $1/4(1-\a)\,\le\, R(\a)\,\le\, \a/\sin^2(\pi \a)$ for all $\a \in [1/2, 1],$ as desired.

\end{proof}

\begin{Vier} {\em (a) As mentioned before one has $G_{1,t}^{(n)} (\lbd) = (-1)^n e^{-\lbd} (\lbd + t -n).$ On the other hand the simple formula 
$(\lbd^\a)^{(n)}\; = \Gamma(\a +1)\lbd^{\a-n}/\Gamma(\a +1-n)$
holds, and one could hence deploy Faa di Bruno's formula in order to investigate the CM property of $\Fat (\lbd) = G_{1,t} (\lbd^\a)$ and try to find the function $R(\a)$ explicitly for all $\a\in [0,1].$ The combinatorial character of the resulting formul\ae\, is however somewhat mysterious, and does not seem to give any simple reason why $R$ is the identity function on $[0,1/2]$ and suddenly behaves differently when $\a > 1/2.$ An exact formula for $R$ on $(1/2,1],$ if any, is probably not given in terms of elementary functions.

\vspace{2mm}

(b) During the proof of the theorem we will establish the continuity of $R,$ which is hence an homeomorphism from $[0,1]$ to $[0,+\infty].$ It is plausible that  $\lim_{\a\to 1} (1-\a)R(\a)$ exists in $[1/4,1]$ and that $R$ is convex.

\vspace{2mm}

(c) In the above proof one sees that the curious equivalence
$$\Fat\;\mbox{is log-convex}\;\Longleftrightarrow\;\Fat\;\mbox{is CM}\;\Longleftrightarrow\; \Hta\;\mbox{is Bernstein}$$ 
holds for all $\a\le 1/2.$ We stress that this is not true anymore when $\a > 1/2.$ For the first equivalence we have for example $G_{3/4, t}$ is log-convex $\Leftrightarrow\, t\ge 4/3$ and on the other hand, although
$$(-1)^n G_{3/4, 4/3}^{(n)} (\lbd) \ge 0$$ 
for all $n\le 4,$ a further computation yields 
$$(-1)^5 G_{3/4,4/3}^{(5)} (\lbd)\; =\; \frac{3 e^{-\lbd^{3/4}}}{4^5\lbd^{1/2}} (195 \mu^5 + 35 \mu^4 - 150 \mu^3 - 135 \mu^2 - 27 \mu +81)$$
with the notation $\mu = \lbd^{-3/4}:$ one sees that the expression on the right-hand side is negative for $\mu = 4/5.$ For the more subtle invalidity of the second equivalence when $\a > 1/2,$ we refer to Proposition 7 (a) below.}
\end{Vier}

\subsection{End of the proof} Part (a) is an easy consequence of Kanter's representation (\ref{Kant}) and Lemma 1. The latter entails indeed clearly  that the function $b_\a^{-r/\a}(u)$ is increasing and convex for all $r >0,$ resp. decreasing and concave for all $r\in [-\a,0).$  In particular,  the variable $b_\a^{-r/\a}(U)$ is monotone for all $r\ge -\a$ and $r\neq 0$,  because its distribution function is concave for all $r >0$ resp. convex for all $r\in [-\a,0)$. On the other hand, one sees explicitly that the function $x\mapsto \lbd_{\a,r} (e^x)$ is log-concave for all $\a\in (0,1)$ and $r\neq 0,$ where $\lbd_{\a,r}$ is the density of $L^{r(\a-1)/\a}.$ By Theorem 3.7. in \cite{CT}, this means that $L^{r(\a-1)/\a}$ is  multiplicatively strong unimodal, in other words, that its independent product with any unimodal random variable remains unimodal. Recalling
$$\Zar\; \elaw\; L^{r(\a-1)/\a}\times b_\a^{-r/\a}(U),$$
we can deduce that $\Zar\in\UU$ as soon as $r \ge -\a.$ 

For the non-monotonicity assertion we simply remark from the Humbert-Pollard representation - see e.g. (14.31) in \cite{S}:
\begin{equation}
\label{HP}
\fa(x)\; =\; \sum_{n\ge 1} \frac{(-1)^{n-1}\Gamma(1+\a n)\sin \pi\a n}{\pi n!}\, x^{-\a n -1}\; =\; \sum_{n\ge 1} \frac{(-1)^n}{n!\Gamma(-\a n)}\, x^{-\a n -1},\quad x>0,
\end{equation}
Linnik's asymptotic $\lim_{t\to 0} x^{\a/(1-\a)} \log \fa(x) = c_\a \in (-\infty,0)$ - see e.g. (14.35) in \cite{S}, and a change of variable, that $\far(0+) = 0$ as soon as $r >-\a,$ whence $\Zar\not\in \MM$ as desired.  

\vspace{2mm}

(b) From above, we already know that $Z_\a^{-\a}\in\UU,$ and it is hence enough to show that zero is a mode of $Z_\a^{-\a}$ for $\a \le 1/2$ and not a mode of $Z_\a^{-\a}$ for $\a >1/2.$ When $\a = 1/2$ this comes from the explicit formula 
$$f_{1/2}^{-1/2}(x)\; = \; \frac{1}{\sqrt{\pi}}e^{-x^2/4}$$ 
and when $\a \neq 1/2$ one obtains the desired properties easily from the limiting values 
$$f_\a^{-\a}(0+)\; =\; \frac{1}{\Gamma(1-\a)}\quad \mbox{and}\quad (f_\a^{-\a})'(0+)\; =\; \frac{-2}{\Gamma(1-2\a)}$$
which were derived during the introduction.

\vspace{2mm}

(c) We first recall the limiting value $\far(0+) = +\infty$ for all $r <-\a$ - see again the end of the introduction, whence the first equivalence 
$$\Zar \in\UU\;\Longleftrightarrow\;\Zar\in\MM.$$ 
Moreover, a change of variable yields
$$r^2(\far)'(x)\; =\; -x^{-(1-1/r)}\har(x^{1/r})$$
with $\har(x) = (1-r) \fa(x) + x \fa'(x).$ From (\ref{HP}) and (14.35) in \cite{S} one sees that $\fa$ and $\har$ are bounded and integrable. Using $\fa(0) = 0,$ an integration by parts yields
\begin{eqnarray*}
\int_{\rl^+} e^{-\lbd x} \har(x) dx 
& = & -r e^{-\lbd^\a}\, +\,\lbd\int_{\rl^+} x e^{-\lbd x} \fa(x) dx\\
& = & -r e^{-\lbd^\a}\, +\,\lbd \,\frac{{\rm d}}{{\rm d}\lbd}\lpa e^{-\lbd^a}\rpa \; =\; (\a \lbd^\a -r)e^{-\lbd^a} \;=\;  \Far(\lbd)
\end{eqnarray*}
for all $\lbd \ge 0.$ Berstein's theorem and Lemma 2 entail then the two other desired equivalences
\begin{equation}
\label{Gew}
\Zar\,\in\,\MM\;\Leftrightarrow\; \har \ge 0\;\Leftrightarrow\;\Far\;\, \mbox{is CM}\;\Leftrightarrow\;\a\le 1/2\;\, {\rm or}\;\, r \le -R(\a),
\end{equation}
with the above notation for $R(\a).$ From $\har\ge 0\Leftrightarrow  r \le -R(\a)$ and the continuity $(\a, r) \mapsto \har (x)$ for all $x>0$ - which is itself a clear consequence of (\ref{HP}) and its derivative - we get the continuity of $\a\mapsto R(\a)$, so that $R$ is homeomorphic from $[0,1]$ onto $[0,+\infty].$ Last, the remaining bound $R(\a)\le \a/(1-\a)$ for $\a \ge 1/2$ comes from $\Zar\,\in\,\MM\,\Leftrightarrow\, r \le -R(\a)$ and the main theorem of \cite{PSS}.
\fin

\begin{Fuenf} {\em (a) This result shows in particular that $\Zar$ is unimodal for all $r\in\rl$ and $\a\le 1/2.$ As a matter of fact, this property is also a direct consequence of our previous theorem in \cite{Th}, which established the equivalence $\Za$ is MSU $\Leftrightarrow \a\le 1/2$ - see p. 2 in \cite{Th} for a definition of the MSU property. From (1.2) in \cite{Th} it is namely clear that the equivalence 
$$\Za\,\mbox{is MSU}\; \Longleftrightarrow\; \Zar\; \mbox{is MSU for all $r\in \rl^*$}$$ 
holds, and Theorem 3.6. in \cite{CT} entails the implication $X$ is MSU $\Rightarrow X\in\UU,$ for all positive random variables $X$. The theorem of \cite{Th} gives hence another proof that $R(\a) =\a$ whenever $\a\le 1/2.$ The argument which is given above in terms of log-convexity and  Bernstein's property is perhaps more illuminating, though also more surprising - recall Remark 4 (c) above. 

\vspace{2mm}

(b) Our two theorems in \cite{Th} and in the present article raise the natural question whether the equivalence
$$X\; \mbox{is MSU}\;\Longleftrightarrow\; X^r\,\in\,\UU\;\;\mbox{for all $r\in\rl$}$$
should not hold for all positive random variables, since it is true for the stable ones. Before studying this conjecture it would be interesting to investigate the MSU property for a larger class of positive random variables than the stable family, for example the self-decomposable subclass. One may ask if this property could not be characterized in terms of the L\'evy measure. 

\vspace{2mm}

(c) The function 
$$x\; \mapsto\; \frac{(\a x^\a +r)e^{-x^\a}}{(1+r)\Gamma(1+1/\a)}$$
is a density on $\rl ^+$ for all $\a,r > 0$ and Lemma 3 entails that it is log-convex for all $r\ge 1/4$ and $\a\le 1- 1/4r.$ By Theorem 51.4 in \cite{S}, it is then the density of an ID law, and this property is also easy to verify for $\a = 1$ and all $r > 0.$ However it is not clear whether the law of the corresponding variable remains ID for all $\a \in ]1-1/4r, 1[$, like for the CM density function $e^{-x^\a}/\Gamma(1+1/\a).$}
\end{Fuenf}

\section{Further results and remarks}

Considering the smooth function $\gar(x) = (1+1/r) \fa(x) + (x/r) \fa'(x)$ for all $\a\in (0,1), r >0,$ it is clear from the above proof of Part (c) of the theorem that the equivalences
$$\gar\;\mbox{is a density on $\rl^+$}\;\Longleftrightarrow\; Z_\a^{-r}\,\in\,\MM\;\Longleftrightarrow\; r\,\ge\, R(\a)$$
hold, and we will set $\Xar$ for the random variable with density $\gar,$ whenever it exists. The following corollary provides a multiplicative factorization somehow analogous to the main theorem of \cite{PSS}. 

\begin{Korse} One has the equivalence
\begin{equation}
\label{MrK}r\,\ge\, R(\a)\;\Longleftrightarrow\; \Zar\,\elaw\,e^L\,\times\, X,
\end{equation} 
where $L\sim$ {\em Exp (1)} and $X$ is an independent random variable. Moreover, when (\ref{MrK}) holds, then $X\elaw X_{\a,r}^r.$
\end{Korse}

\begin{proof} From the identity $e^L\elaw U^{-1}$ with $U\sim$ Unif(0,1), the required equivalence reads
$$r\,\le\, -R(\a)\;\Longleftrightarrow\; \Zar\,\elaw\, U\,\times\, Y$$
with $Y\perp U.$ The latter is however a direct consequence of Khintchine's theorem - see e.g. Footnote 22 p. 155 in \cite{F} - since we have previously shown $r\,\le\, -R(\a)\;\Leftrightarrow\; \Zar\,\in\,\UU_0.$
For the remaining identity in law we first remark from (\ref{MrK}) that for all $s < \a/r\le 1$ one has 
\begin{eqnarray*}
\EE[X^s] \; = \; (1-s) \EE[\Za^{rs}] & = & 2 \EE[\Za^{rs}] \; -\; (1+s)  \EE[\Za^{rs}]\\
& = & \int_0^\infty \!\!\! x^s(2\far(x) + x(\far)'(x)) dx\\
& = & \int_0^\infty\!\!\! x^s (x^{1/r-1} \gar(x^{1/r})/r) dx\; =\;\EE[X_{\a,r}^{rs}],
\end{eqnarray*}
where the third equation comes from an integration by parts - noticing that both boundary values vanish from $r > 0$ and Linnik's asymptotic  resp. from $rs <\a$ and the Humbert-Pollard representation - and the fourth one from a direct computation. Mellin inversion entails the desired identity 

$$X\;\elaw\; X_{\a,r}^r.$$

\end{proof}

The Laplace transform of $\Xar$ reads 
$$\EE[e^{-\lbd\Xar}]\; =\; (1+\a\lbd^\a/r)e^{-\lbd^\a}$$
and because $r\ge \a,$ we know from the proof of Lemma 3 that $\Xar$ is ID as soon as $\a\le 1/2.$ It is interesting to ask whether its is always the case, in other words, whether the second equivalence in Remark 4 (c) always holds. From Corollary 6 and the main theorem of \cite{PSS} one obtains namely
$$ e^{L/r}\;\times\; \Xar\;\elaw\;\Za\;\elaw\; L^{-1/r}\;\times\; e^{\Yar}$$
as soon as $r\ge \a/(1-\a)\ge R(\a),$ with $\Yar$ always an ID random variable. If $\Xar$ is also ID, one obtains through $\Za$ a kind of cross-correpondence between Exp and some ID laws which is perhaps more than a formal one. Besides, the infinite divisibility of $X_{\a, R(\a)}$ for $\a >1/2$ would give a characterization of $\Zar\in\MM$ in terms of the Bernstein property for $-\log \Far,$ would be more tractable than the CM property for $\Far$ - see Remark 8 (a) below - in investigating further properties of the frontier function $R(\a)$. The following proposition shows however that $\Xar$ is not always ID, at least for $\a > 1/2$ and $r =R(\a).$ For completeness we also discuss the self-decomposability of $\Xar.$ 

\begin{Prosi} {\em (a)}  There exists a homeomorphism $\Rde : [0,1]\to [0,+\infty]$ such that
$$\Xar\;\mbox{is ID}\;\Longleftrightarrow\; r\ge \Rde(\a).$$
Moreover one has $\Rde(\a) = \a$ for all $\a\le 1/2$ and $\a/\sin^2(\pi\a)\ge\Rde(\a) > R(\a)$ for all $\a > 1/2.$

\vspace{1mm}

{\em (b)} There exists a homeomorphism $\Rta : [0,1]\to [0,+\infty]$ such that
$$\Xar\;\mbox{is SD}\;\Longleftrightarrow\; r\ge \Rta(\a).$$
Moreover one has $\Rta(\a) = \a$ for all $\a\le 1/2$ and $\a/\sin^2(\pi\a)\ge\Rta(\a) > \Rde(\a)$ for all $\a > 1/2.$
\end{Prosi}

\begin{proof} (a) The case $\a\le 1/2$ comes directly from Lemma 3 and Remark 4 (c). In the case $\a > 1/2$ one sees from the proof of Lemma 3 that the equivalences 
$$\Xar\;\mbox{is ID}\;\Longleftrightarrow\;\Hta\;\mbox{is Bernstein}\;\Longleftrightarrow\; \sup\{U_\a (x), \, x\ge 0\}\; \le \; t$$
hold, with the notation of Lemma 3 and $t = r/\a.$ Setting
$$\Tde(\a) =\inf\{t\ge 0\;/\; \Hta\;\mbox{is Bernstein}\},$$ the above second equivalence and the smoothness of $(\a, x) \mapsto \Ua(x)$ entail the continuity of $\a\mapsto\Tde(\a).$ It is not clear from the definition of $\Ua$ that the function $\Tde$ does not decrease on $[0,1],$ but instead one can prove this indirectly  just like in Lemma 3, since $\lbd\mapsto\psi(\lbd^\b)$ is Bernstein for all Bernstein functions $\psi$ and $0<\b < 1.$ The function $\Rde(\a) = \a\Tde (\a)\ge R(\a)$ yields the required homeomorphism and the upper bound $\Rde(\a)\le \a/\sin^2(\pi\a)$ comes from the integral representation of $\Ua$ in Lemma 3. 

We finally show $\Rde(\a) > R(\a),$ in other words, that $\Xda =X_{\a, R(\a)}$ is not ID as soon as $\a > 1/2.$ Fix $\a\in (1/2,1)$ and set $\ga  = g_\a^{R(\a)}$ for the density of $\Xda.$ From (\ref{HP}) we first get
$$\gar (x) = \sum_{n\ge 1} \frac{(-1)^n (1- \a n/r)}{n!\Gamma(-\a n)} \,x^{-\a n -1},\quad x,r > 0,$$
whence we easily deduce the existence of $x_0 > 0$ and of an open neighbourhood $\cV$ of $R(\a)>\a$ such that $\gar(x) > 0$ for all $r \in\cV$ and all $x > x_0.$ On the other hand, from the positivity of $\fa$ and $\fa'$ on $(0, m_\a]$ - here $m_\a$ denotes the mode of $\fa$, there exists $x_1 >m_\a > 0$ such that $\gar(x) > 0$ for all $r > 0$ and $x\in (0, x_1).$ By (\ref{Gew}) however there must exist, for all $r < R(\a)$ close enough to $R(\a),$ some $x_r \in [x_1, x_0]$ such that $\gar(x_r) < 0.$ By a straightforward compacity argument, this entails that $\ga$ must vanish on $[x_1, x_0].$

Now if $\Xda$ were ID, then its L\'evy-Khintchine exponent would be $H_\a^{T(\a)}$ and one sees from (\ref{MrL}) that $\Xda$ would be driftless and that the support of its infinite (because $\Ua(x) \to 0$ as $x\to 0$) L\'evy measure would contain zero. From Theorem 24.10 in \cite{S} we deduce that the support of $\Xda$ would be the whole $\rl^+,$ and Theorem 3.3. in \cite{Sh} entails that $\ga$ would never vanish on $]0,+\infty[,$ a contradiction. 

\vspace{2mm}

(b) From the proof of Lemma 3 and Corollary 15.11 in \cite{S} we see that $\Xar$ is SD if and only if the function
$$x\;\mapsto\; x^{-\a} - \Gamma(1-\a) \Ea (-x^\a t)$$ 
is non-increasing, with $t = r/\a,$ which is equivalent to the upper bound $\Va(x)\le t$ for all $x \ge 0,$ where we have set $\Va(x) = x^{2\a} \Gamma(1-\a) \Ea'(-x^\a).$ As above we can hence introduce the function
$$\Rta(\a) =\inf\{ r > 0 \; /\; \mbox{$\Xar$ is SD}\},$$ 
whose continuity on $[0,1]$ is clear, and which yields the desired characterization. One shows that $\Tta (\a) =\Rta(\a)/\a$ does not decrease  indirectly like in Part (a), using this time Proposition 4.1. in \cite{SSS}, so that $\Rta$ increases and is a homeomorphism. Differentiating (\ref{ML}) we finally obtain
$$\Va(x)\; =\; \frac{1}{\Gamma(\a +1)}\int_0^\infty \!\frac{u^\a e^{-u} }{(u/x)^{2\a} + 2(u/x)^\a\cos \pi\a +1}\, du,$$
whence $\Tta(\a) = 1\,\Leftrightarrow\,\Rta(\a) =\a$ for all $\a \le 1/2,$ respectively $\Tta(\a) \le 1/\sin^2(\pi\a)\,\Leftrightarrow\,\Rta(\a) \le\a/\sin^2(\pi\a)$ for all $\a > 1/2.$ For the strict lower bound in the case $\a >1/2$ we first recall from (a) that 
$$\Rta(\a) \;=\;\a\sup\{\Ua(x), \; x\ge 0\}\; >\; \a.$$ 
Since $\Ua(0) = 0$ and $\Ua(x) \to 1$ as $x\to \infty,$ the corresponding maximal value is attained in $]0, +\infty[$ and the function $$x^{-\a}(\Tde(\a) - \Ua (\Tde(\a)x^\a))\; =\; x^{-\a} \, -\, \Gamma(1-\a) \Ea(-x^\a \Tde(\a))$$ 
must vanish inside $]0, +\infty[$. Hence the latter function must increase somewhere in $]0, +\infty[,$ since it is clearly positive in the neighbourhood of $+\infty.$ This shows that the ID variable $X_{\a, \Rde(\a)}$ is not SD and one has $\Rta(\a) > \Rde(\a)$ for all $\a > 1/2.$

\end{proof}

\begin{Acht} {\em (a) As mentioned before, the limiting functions $\Rde$ and $\Rta$ can be expressed as extremal values of an explicit function: one has
$$\Rde(\a)\; =\; \a\sup\{\Ua(x), \; x\ge 0\}\quad\mbox{and}\quad \Rta(\a)\; =\; \a\sup\{\Va(x), \; x\ge 0\},$$
and the classical integral representations for $\Ea$ and $\Ea'$ entail readily $\Rde(\a) =\Rta(\a) = \a$ for all $\a\le 1/2.$ As for $R(\a),$ these representations give however little hope of a closed formula for $\Rde(\a)$ and $\Rta(\a)$ when $\a >1/2$.

\vspace{2mm}

(b) When $r\ge R(\a)$ writing $e^{-\lbd^\a} = (1+ \a\lbd^\a/r)^{-1} (F_\a^{-r}(\lbd)/r)$ provides another, additive factorization for $\Za$:
$$\Za\; \elaw\; (\a/r)^{1/\a}\Ma\; +\; \Xar,$$
where $\Ma$ is a so-called Mittag-Leffler variable - see \cite{P} - and $\Ma\perp\Xar.$ Notice that a direct moment computation - see also the final remark in \cite{P} - gives the factorization
$$\Ma\;\elaw\; \Za\;\times\; L^{1/\a}$$ 
for all $\a\in (0,1)$, with $L\sim$ Exp(1) and $L\perp\Ma.$ It is however purposeless to search for a reverse multiplicative factorisation of $\Za$ through $\Ma$ because $\Ma\in\UU_0$ (its density function is namely $\a x^{\a -1} \Ea'(-x^\a)$ which is non-increasing on $\rl^+$ - see (18.1.6) in \cite{E}) and $\Za\not\in\UU_0$, which would contradict Khintchine's theorem.}
\end{Acht}

\bigskip

\noindent
{\bf Aknowledgement.} This paper was supported by the grant ANR-09-BLAN-0084-01.

\end{document}